\documentclass[letterpaper, 10 pt, conference]{ieeeconf}  

\IEEEoverridecommandlockouts                              
\usepackage{graphicx} 
\usepackage{amsmath} 
\usepackage{amssymb}  

\usepackage{algorithmic}
\usepackage{array}
\usepackage{color}
\usepackage{subfig}

\newcommand{\real}{{\mathbb{R}}}
\newcommand{\subj}{\text{subj. to}}

\title{\LARGE \bf
Computing minimum-time trajectories for quadrotors \\ via transverse coordinates}

\author{Sara~Spedicato
	and~Giuseppe~Notarstefano$^{*}$
	\thanks{$^{*}$ S. Spedicato and G. Notarstefano are with the Department of Engineering, Universit\`a del Salento, Via per Monteroni, 73100 Lecce, Italy, e-mail: \tt\small{name.lastname@unisalento.it.}}%
	\thanks{This result is part of a project that has received funding from the European Research Council (ERC) under the European Union's Horizon 2020 research and innovation programme (grant agreement No 638992 - OPT4SMART).}}

\begin{document}

\maketitle
\thispagestyle{empty}
\pagestyle{empty}

\begin{abstract}
  In this paper we present a novel strategy to compute minimum-time trajectories
  for quadrotors. In particular, we consider the motion in constrained
  environments, taking into account the physical limitations of the
  vehicle. Instead of approaching the optimization problem in its standard
  time-parameterized formulation, the proposed strategy is based on an appealing
  re-formulation. Transverse coordinates, expressing the distance from a ``reference" path, are used to parameterize the vehicle position and a
  spatial parameter is used as independent variable. This re-formulation allows
  us to (i) obtain a fixed horizon problem and (ii) easily formulate (even
  complex) position constraints. The effectiveness of the proposed
  strategy is proven by numerical computations on two different illustrative
  scenarios.
\end{abstract}

\section{INTRODUCTION}
Numerous applications involving quadrotors
require them to move inside areas characterized by
physical boundaries, obstacles and even tight space constraints (as e.g., urban
environments) in order to accomplish their robotics tasks.  
Trajectory generation, a core step for physical
task realization \cite{dadkhah2012survey}, becomes extremely challenging in this
scenario.  A physically realizable trajectory must satisfy (i) the (nonlinear)
system dynamics, (ii) the physical limits of the vehicle, 
and (iii) the position constraints. 
Moreover, the computation of optimal (rather than just feasible) trajectories is more and more becoming a
necessity to ensure an efficient accomplishment of the tasks, posing an additional challenge in the trajectory generation problem.
 
Quadrotor trajectory planning in constrained environments can be addressed using
a reacting (e.g.,  \cite{hou2016dynamic, furci2015plan}) or a planning approach.
The majority of the planning algorithms, such as
\cite{mellinger2011minimum},
\cite{bouktir2008trajectory}, \cite{van2013time}, take advantage of the
differential flatness property to avoid the integration of nonlinear
differential equations, reduce the order of the problem and simplify the
definition of position constraints.
The problem is posed in the flat output space, where outputs are approximated as a combination of primitives, 
such as polynomial functions
\cite{mellinger2011minimum}, B-splines
\cite{bouktir2008trajectory}, or convex functions \cite{van2013time}. 
The parameters of the primitives are the optimization variables.
When dealing with obstacle dense environments, trajectory generation is often
performed using a decoupled approach (\cite{bry2015aggressive},
\cite{koyuncu2008probabilistic}, \cite{bouffard2009hybrid}) aiming to avoid the
difficult definition of spatial constraints and the high computational cost.  In
a first stage, a collision-free path is generated by sampling-based path planning
algorithms, such as the Rapidly-exploring Random Tree (RRT) in
\cite{bry2015aggressive, koyuncu2008probabilistic} or the Probabilistic Roadmap
(PRM) in \cite{bouffard2009hybrid}, and without the dynamics constraint.  In a
second stage, an optimal trajectory (satisfying the system dynamics) is
generated from the collision free path. Optimization techniques such as
\cite{cowling2010direct}, \cite{mellinger2011minimum},
\cite{bouktir2008trajectory} can be used at this stage. 
Differently from these works, the one in \cite{van2013time}
considers a space-parameterized problem re-formulation, suitable for modeling
complex flight scenarios, and thus avoids the decoupled approach.
Finally, only few works do not consider the
differential flatness hypothesis, such as
\cite{augugliaro2012generation}.
Nevertheless, in \cite{augugliaro2012generation} a discretized simple point-mass
dynamics and approximated convex constraints are considered.
	
Our main contribution is the design of an optimization framework to generate
feasible minimum-time quadrotor trajectories in geometrically constrained
environments.  In our problem setup (i) the vehicle nonlinear dynamics is
directly included into the optimization problem (the differential flatness
assumption is not required) and (ii) the vehicle position, which is an
optimization variable, is not approximated via primitives.
Similarly to \cite{van2013time}, the minimum-time
problem is formulated by using the path parameter as independent variable: a fixed
horizon problem is obtained and even complex position
constraints can be easily formulated.
The optimization problem is solved combining the Projection Operator Newton
method for Trajectory Optimization (PRONTO) \cite{hauser2002projection} with a
barrier function approach \cite{hauser2006barrier}.
Similar problem
setups and optimization methods are presented in \cite{hauslerenergy} for
two-wheeled mobile robots and in \cite{rucco2015virtual} for UAVs. In
\cite{hauslerenergy} a minimum-energy problem in a bi-dimensional space with
obstacles is addressed by considering a standard formulation (i.e., with
variables parameterized by the time).
In \cite{rucco2015virtual} a virtual target approach (without position constraints) is presented using an
arc-length parameterized transverse dynamics model. 
Differently from these works, we solve
the optimization problem in the three-dimensional space and we define obstacle
constraints by an arc-length parameterization.  As an additional contribution,
we present numerical computations on two challenging scenarios with respectively, a circular tube of variable section and a narrow corridor with an obstacle, as constrained spaces.

The paper is organized as follows. In Section \ref{sec:problem_formulation} we
present the standard formulation of the optimization problem we aim to solve. In
Section \ref{sec:strategy} our trajectory generation strategy, based on an
appealing reformulation of the problem, is described.  Finally, in Section
\ref{sec:computations}, we provide numerical computations discussing the
minimum-time vehicle behavior.
	
\section{The quadrotor minimum-time problem}
\label{sec:problem_formulation}
We first briefly introduce the quadrotor model used in the paper and then recall
the standard problem formulation.
		
The quadrotor dynamics can be described by the so called vectored-thrust
dynamical model, \cite{hua2013introduction},
\begin{align}
  \dot{\pmb{p}} &= \pmb{\text{v}}\label{eq:state_pos}\\
  \pmb{\dot{\text{v}}}& = g \pmb{e}_3 - fm^{-1} R(\pmb{\Phi})  \pmb{e}_3 \label{eq:state_vel}\\
  \dot{\pmb{\Phi}}&=J \pmb{\omega} \label{eq:state_ang}\\
  \pmb{\dot{\omega}} &= -I^{-1}\Omega I \pmb{\omega} +
                       I^{-1} \pmb{\gamma}, \label{eq:state_omega}
\end{align}
where $\pmb{p}=[p_1 \; p_2 \: p_3]^T$ is the vehicle position expressed in the
inertial frame $\mathcal{F}_i$,
$\pmb{\text{v}}=[\text{v}_1 \; \text{v}_{2} \; \text{v}_3]^T$ is the linear
velocity expressed in $\mathcal{F}_i$, $\pmb{\omega} = [\omega_1 \; \omega_2 \; \omega_3]^T$ is the
angular rate expressed in the body frame $\mathcal{F}_b$ ($\Omega \in so(3)$ is
such that, for
$\pmb{\beta} \in \mathbb{R}^3, \Omega \pmb{\beta} = \pmb{\omega} \times
\pmb{\beta}$),
and $J \in \mathbb{R}^{3 \times 3}$ is the Jacobian matrix.  The rotation matrix
$R(\pmb{\Phi}) \in SO(3)$, mapping vectors in $\mathcal{F}_b$ into vectors in
$\mathcal{F}_i$, is parameterized using roll-pitch-yaw angles, denoted by
$\pmb{\Phi} := [\varphi \; \theta \; \psi]^T$.
Furthermore, $m$ and $I$ are respectively the mass and the inertia matrix, $g$
is the gravity constant, and $\pmb{e}_3=[ 0 \; 0 \; 1]^T$. The vehicle is
controlled by the thrust $f$ and the torques $\pmb{\gamma}$.
	
For the vehicle maneuvering, we use a cascade control scheme with an off-board
position/attitude control loop and an on-board angular rate controller.
Assuming that the \emph{virtual} control input $\pmb{\omega}$ is tracked by the
on-board angular rate controller, we restrict our trajectory generation problem
on the position/attitude subsystem (\ref{eq:state_pos}-\ref{eq:state_ang}),
which can be written in state-space form as
\begin{equation}
  \dot{\pmb{x}}(t) = f(\pmb{x}(t),\pmb{u}(t)),
  \label{eq:state_space}
\end{equation}
with state $\pmb{x} = [\pmb{p}^T \; \pmb{\text{v}}^T \; \pmb{\Phi}^T]^T$ and
input $\pmb{u} = [\pmb{\omega}^T \; f]^T$.
	
We deal with the following optimal control problem:
  \begin{align}
    \min_{\pmb{x}(\cdot),\pmb{u}(\cdot)} &\;  \int_0^T \!\!\!1 \; d\tau \nonumber \\
    \!\!\subj &\; \dot{\pmb{x}}(t) = f(\pmb{x}(t),\pmb{u}(t)), \quad \pmb{x}(0)
    = \pmb{x_0} \quad \text{\emph{(dynamics)}} \nonumber\\
    & \; |\omega_i(t)| \leq \omega_{i,max}, \; i = 1,2,3, \; \text{\emph{(angular rate)}} \nonumber \\
    & \;  0 < f_{min} \leq f(t) \leq f_{max} \; \text{\emph{(thrust)}} \nonumber \\
    & \; |\varphi(t)| \leq \varphi_{max}, \; |\theta(t)| \leq \theta_{max} \; \text{\emph{(angles)}}\nonumber \\
   & \; c_{obs}(\pmb{p}(t)) \leq 0 \; \text{\emph{(position constraints)}} \label{eq:mintime_standard}
\end{align}
where $\omega_{i,max}$ are the bounds on the angular rate, $f_{min}$ and $f_{max}$ are lower an upper bounds on thrust,
$\varphi_{max}$ and $\theta_{max}$ are the bounds on the roll and pitch angles respectivelly, and
$c_{obs} : \real^3 \rightarrow \real$ represents position constraints taking into account physical boundaries and obstacles.
The position constraints may also represent GPS denied areas or spaces with limited communication.
Angular rates have bounds, due to, for example, a limited range of the
gyroscopes or to controller limitations. The vehicle thrust is also limited:
quadrotor vehicles can only generate positive thrust and the maximum rotor speed
is limited. Furthermore, constraints on angles are imposed into the optimization
problem in order to avoid acrobatic vehicle configurations.
	
\section{Minimum-Time Trajectory Generation Strategy}
\label{sec:strategy}
	
In order to solve the minimum-time problem, we consider an equivalent
formulation of \eqref{eq:mintime_standard} where a spatial parameter (instead of time) is the independent
variable. 
The derivation of the equivalent problem formulation is based on three steps:
(i) generation of a ``reference'' path, (ii) definition of the transverse
dynamics, in which the vehicle position is parameterized by transverse coordinates (with respect to that path), and (iii) mapping of
time-parameterized cost and constraints into arc-length parameterized ones.  The
PRONTO method, combined with a barrier function approach, is
then used to numerically solve the optimization problem.
\subsection{Minimum-Time problem: equivalent formulation}
The first step to convert the standard problem into the new formulation is the computation of a ``reference'' path $\bar{\boldsymbol{p}}_{r}(s)$, $\forall s \in [0,L]$, where
$s$ is the arc-length of the path and $L$ is its total length.  In the
following, we denote the arc-length parameterized functions with a bar, and the
derivative with respect to the arc-length with a prime, i.e.,
$\bar{\boldsymbol{p}}_{r}'(s) := d\bar{\boldsymbol{p}}_{r}(s)/ds$.
The ``reference'' path is computed as a $C^{\infty}$ geometric curve, e.g., using
arctangent functions, as in our numerical computations. 
Note that $\bar{\pmb{p}}_{r}(\cdot)$ is not required to 
satisfy the position
constraints, and it can be suitably chosen.

With $\bar{\pmb{p}}_{r}(\cdot)$ in hand, 
we re-write the system dynamics \eqref{eq:state_space} as follows.

\begin{figure}[htbp]
  \begin{center}
    \includegraphics[scale=1.3]{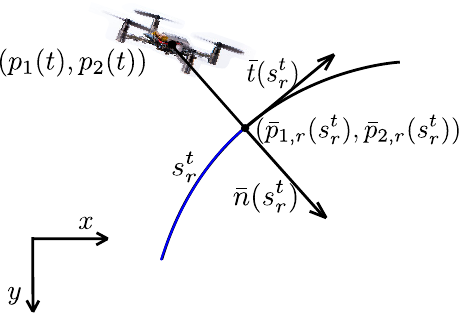}
      \caption{Selection of the arc-length $s_{r}^t$ identifying the point on the ``reference'' path at minimum distance from the quadrotor position at the time instant t.}
    \label{fig:manreg}
  \end{center}
\end{figure}

First, we design a change of coordinates from the inertial position $\boldsymbol{p} \in \real^3$ to
the transverse coordinate vector $\boldsymbol{w} \in \real^2$, such that $\boldsymbol{w} = [w_1 \; w_2]^T$. Let us consider the ``reference'' path $\bar{\boldsymbol{p}}_r(\cdot)$.
A Serret-Frenet frame, whose origin has $\bar{\boldsymbol{p}}_{r}(s)$ as coordinates, can be defined $\forall s \in [0,L]$. 
In particular, the tangent, normal and bi-normal vectors, respectively $\bar{\boldsymbol{t}}(s), \bar{\boldsymbol{n}}(s),
\bar{\boldsymbol{b}}(s)$, are defined, with components in the inertial frame, as
\begin{align}
  \bar{\boldsymbol{t}}(s)&:=\bar{\boldsymbol{p}}'_{r}(s),\\
  \bar{\boldsymbol{n}}(s)&:=\frac{\bar{\boldsymbol{p}}''_{r}(s)}{\bar{k}(s)},\\
  \bar{\boldsymbol{b}}(s)&:= \bar{\boldsymbol{t}}(s) \times
                           \bar{\boldsymbol{n}}(s),
\end{align}
where $\bar{k}(s):=\Vert \bar{\boldsymbol{p}}''_{r}(s)\Vert_2$ is the
curvature of $\bar{\boldsymbol{p}}_{r}(\cdot)$ at $s$.
Moreover, we define the rotation matrix $\bar{R}_{SF}:=[\bar{\boldsymbol{t}} \; \bar{\boldsymbol{n}} \; \bar{\boldsymbol{b}}]$ mapping
vectors with components in the Serret-Frenet frame into vectors with
components in the inertial frame. 
According to the Serret-Frenet formulas, the arc-length derivative of the Serret-Frenet rotation matrix is 
\begin{equation}
  \bar{R}'_{SF}(s) =
  \bar{R}_{SF}(s)
  \left[
    \begin{array}{ccc}
      0 & -\bar{k}(s)  & 0 \\
      \bar{k}(s) & 0 & -\bar{\tau}(s) \\
      0 & \bar{\tau}(s) & 0 \\
    \end{array}
  \right],
  \label{eq:Rsf'*dot_s}
\end{equation}
where $\bar{\tau}(s):=\bar{\boldsymbol{n}}(s) \; \bar{\boldsymbol{b}}'(s)$ is
the torsion of $\bar{\boldsymbol{p}}_{r}(\cdot)$ at $s$.
Let us now consider the quadrotor center of mass with position $\boldsymbol{p}$.
As depicted in Figure \ref{fig:manreg}, its orthogonal projection on the ``reference'' path identifies a point with position $\bar{\boldsymbol{p}}_r(\pi(\boldsymbol{p}))$,
where the function $\pi : \real_{0}^+ \rightarrow \real_{0}^+$ is defined as
\begin{align}
\pi(\boldsymbol{p}):=\text{arg} \min_{s \in \real_{0}^+}
\|\boldsymbol{p}-\bar{\boldsymbol{p}}_{r}(s)\|^2.
\end{align}
Note that, the minimizing arc-length is unique provided that $\bar{\boldsymbol{p}}_{r}(\cdot)$
is locally a non-intersecting $C^2$ curve with non-vanishing $\bar{\boldsymbol{p}}'_{r} (\cdot)$.
By mapping 
$\boldsymbol{p}-\bar{\boldsymbol{p}}_{r}(\pi(\boldsymbol{p}))$ into a vector with components in the Serret-Frenet frame attached to $\bar{\boldsymbol{p}}_{r}(\pi(\boldsymbol{p}))$, we obtain $\boldsymbol{d}_{SF} := R_{SF}(\pi(\boldsymbol{p}))^T (\boldsymbol{p}-\bar{\boldsymbol{p}}_{r}(\pi(\boldsymbol{p})))$.  
Moreover, noticing that the component related to the tangent vector is always zero by construction, we define the transverse coordinates 
as the second and third components of $\boldsymbol{d}_{SF}$. Thus  $\boldsymbol{d}_{SF} := [0 \; w_1 \; w_2]^T$.
More formally, we define the locally invertible function
$\phi : \real^3 \rightarrow \real^2$, such that $\boldsymbol{w} =  \phi(\boldsymbol{p})$, with components defined as
\begin{align}
  \phi_1(\boldsymbol{p})& := \bar{\boldsymbol{n}}(\pi(\boldsymbol{p}))^T (\boldsymbol{p}-\bar{\boldsymbol{p}}_{r}(\pi(\boldsymbol{p}))),\\
  \phi_2(\boldsymbol{p})& := \bar{\boldsymbol{b}}(\pi(\boldsymbol{p}))^T
                          (\boldsymbol{p}-\bar{\boldsymbol{p}}_{r}(\pi(\boldsymbol{p}))).
                          \label{eq:x_change_coordinates}
\end{align}
Moreover, we denote by $\pmb{\text{v}}_{SF} := \bar{R}_{SF}(\pi(\boldsymbol{p}))^T \pmb{\text{v}}$ the linear velocity of the quadrotor with components in the Serret-Frenet frame. In particular, $\pmb{\text{v}}_{SF} = [\text{v}_t \; \text{v}_n \; \text{v}_b]^T$.

Second, we note that, considering the position function $\boldsymbol{p}(\cdot)$ with the time $t$ as independent variable,
the function $s_r := \pi \circ \boldsymbol{p}$ provides a change of variables from the time $t$ to the arc-length $s$. 
In particular, we denote by $s_r^t := s_r(t)$ the function $s_r(\cdot)$ evaluated at the time $t$, and by $\dot{s}_r^t := \dot{s}_r(t)$ the time derivative of $s_r(\cdot)$, evaluated at $t$. 
A generic arc-length function $\bar{\alpha}(\cdot)$ can thus be expressed as the time function $\bar{\alpha}(s_r(\cdot))$. Note that, the time derivative of $\bar{\alpha}(s_r(\cdot))$
is $\frac{d\bar{\alpha}(s_r(t))}{dt} = \bar{\alpha}'(s_r^t) \; \dot{s}_r^t$. 
Using the function $s_r : \real_{0}^+ \mapsto \real_{0}^+$, the position of the quadrotor center of mass $\boldsymbol{p}(t)$, at time
instant $t$, can be written as
\begin{align}
  \boldsymbol{p}(t) = \bar{\boldsymbol{p}}_{r}(s_{r}^t)+\bar{R}_{SF}(s_{r}^t) \; \boldsymbol{d}_{SF}(t).
  \label{eq:y=y_xi+Rsf*d}
\end{align}
Differentiating \eqref{eq:y=y_xi+Rsf*d} with respect to time, we get
\begin{equation}
  \pmb{\text{v}}(t)=\bar{\boldsymbol{p}}'_{r}(s_{r}^t) \; \dot{s}_{r}^t +
  {\bar{R}}_{SF}'(s_{r}^t) \; \dot{s}_{r}^t \; \boldsymbol{d}_{SF}(t) +
  \bar{R}_{SF}(s_{r}^t) \; \boldsymbol{\dot{d}}_{SF}(t).
  \label{eq:dot_y=y_xi+Rsf*d}
\end{equation}
Multiplying both members of equation \eqref{eq:dot_y=y_xi+Rsf*d} by $\bar{R}^{T}_{SF}$, 
using \eqref{eq:Rsf'*dot_s}, the definition of $\pmb{\text{v}}_{SF}$, and $\bar{\boldsymbol{p}}'_{r}(s_{r}^t) = \bar{R}_{SF}(s_{r}^t) \; 
[1 \; 0 \; 0]^T$,
we get 
\begin{align}
  \dot{s}_{r}^t &= \frac{\text{v}_t(t)}{1-\bar{k}(s_{r}^t) w_1(t)} \label{eq:s_}\\
  \dot{w}_1(t) &= \text{v}_n(t) +\bar{\tau}(s_{r}^t) \dot{s}_{r}^t w_2(t) \label{eq:w1_}\\
  \dot{w}_2(t) &= \text{v}_b(t) -\bar{\tau}(s_{r}^t) \dot{s}_{r}^t w_1(t). \label{eq:w2_}
  \end{align}
Furthermore, since we are using transverse coordinates to parameterize the
vehicle position, we re-write equation \eqref{eq:state_vel} using
$\pmb{\text{v}}_{SF}$ instead of $\pmb{\text{v}}$.
Differentiating both members of
$\pmb{\text{v}}(t) = \bar{R}_{SF}(s_{r}^t) \pmb{\text{v}}_{SF}(t)$ with respect to time,
and using \eqref{eq:state_vel} and \eqref{eq:Rsf'*dot_s}, we get
\begin{align}
  \dot{\pmb{\text{v}}}_{SF} = \bar{R}_{SF}^T \dot{\pmb{\text{v}}} -
  \left[
  \begin{array}{ccc}
    0 & -\bar{k}  & 0 \\
    \bar{k} & 0 & -\bar{\tau} \\
    0 & \bar{\tau} & 0 \\
  \end{array}
  \right] \dot{s}_{r} \pmb{\text{v}}_{SF},
  \label{eq:dot_v_sf}
\end{align}
with $\dot{\pmb{\text{v}}}= ( g \pmb{e}_3 - fm^{-1} R(\pmb{\Phi}) \pmb{e}_3)$.
The dependence on $t$ and $s_r^t$ is omitted for simplicity.
Equations \eqref{eq:w1_}, \eqref{eq:w2_}, \eqref{eq:dot_v_sf} and \eqref{eq:state_ang}
can be compactly written as
\begin{equation}
  \dot{\pmb{q}}_S(t) = \tilde{f}(\pmb{q}_S(t),\pmb{u}(t)),
  \label{eq:dot_qs}
\end{equation}
with state $\pmb{q}_S = [\pmb{w}^T \; \pmb{\text{v}}_{SF}^T \; \pmb{\Phi}^T]^T$
and input $\pmb{u} = [\pmb{\omega}^T \; f]^T$.

Third and final, we write the dynamics \eqref{eq:dot_qs} in terms of the arc-length $s$.
By using the chain rule, we have $\dot{\pmb{q}}_S(t)=\bar{\pmb{q}}_S'(s_{r}^t) \dot{s}_{r}^t$. 
Moreover, we parameterize the state $\pmb{q}_S$ and the input $\boldsymbol{u}$ using
$s_{r}^t$, i.e., $\pmb{q}_S(t)=\bar{\pmb{q}}_S(s_{r}^t)$ and
$\pmb{u}(t)=\bar{\pmb{u}}(s_{r}^t)$. 
Thus, equations
\eqref{eq:dot_qs} can be written as the \emph{transverse dynamics}
\begin{equation*}
  \bar{\pmb{q}}'_S (s) = \bar{f}(\bar{\boldsymbol{q}}_S(s),\bar{\boldsymbol{u}}(s)),
  \label{eq:transv_dynamics}
\end{equation*}
with the arc-length $s$ as independent variable and where
\begin{equation*}
  \bar{f}(\bar{\boldsymbol{q}}_S(s),\bar{\boldsymbol{u}}(s))
  :=
  \tilde{f}(\bar{\boldsymbol{q}}_S(s),\bar{\boldsymbol{u}}(s))
  \frac{1-\bar{k}(s) \bar{w}_1(s)}{\bar{\text{v}}_t(s)}.
\end{equation*}
		
Re-writing the time-parameterized cost and constraints of 
\eqref{eq:mintime_standard}
into arc-length parameterized ones,
the new formulation of the minimum-time problem, using the new set
of coordinates $(\bar{\pmb{q}}_S, \bar{\pmb{u}})$ and with $s$ as independent variable, is
\begin{equation}
  \begin{split}
   \label{eq:mintime}
    \min_{\bar{\pmb{q}}_S(\cdot),\bar{\pmb{u}}(\cdot)} &\;  \int_0^L \!\!\! \quad \frac{1-\bar{k}(s) \bar{w}_1(s)}{\bar{\text{v}}_t(s)} \; ds \\\\
    \!\!\subj &\;
    \bar{\pmb{q}}'_S = \bar{f}(\bar{\boldsymbol{q}}_S(s),\bar{\boldsymbol{u}}(s)), \quad \bar{\pmb{q}}_S(0) = \pmb{q}_{S0} \; \text{\emph{(dynamics)}}\\
    & \; \Big( \frac{\bar{\omega}_i(s)}{\omega_{i, max}} \Big)^2 -1 \leq 0, \; i = 1,2,3
    \; \text{\emph{(angular rate)}}\\
    & \; \Big( \frac{2 \bar{f}(s) - (f_{max}+f_{min})}{(f_{max}-f_{min})}
    \Big)^2 -1 \leq 0
    \; \text{\emph{(thrust)}} \\
    & \; \Big( \frac{\bar{\varphi}(s)}{\bar{\varphi}_{max}(s)} \Big)^2 -1 \leq 0
    \; \text{\emph{(roll angle)}}\\
    & \; \Big( \frac{\bar{\theta}(s)}{\bar{\theta}_{max}(s)} \Big)^2 -1 \leq 0
    \; \text{\emph{(pitch angle)}}\\
    & \; \bar{c}_{obs}(\bar{\boldsymbol{w}}(s)) \leq 0  \;
    \text{\emph{(position constraints).}}
  \end{split}
\end{equation}
The constraints on the roll and pitch angles, the thrust and the angular rates are re-written by means of an $s$-parameterized equivalent form (with smooth functions).
The position constraints $ \bar{c}_{obs}(\bar{\boldsymbol{w}}(s)) \leq 0$ have specific formulations according to the shape of the collision free region.
Considering regions with rectangular or circular
transversal sections, we get the following inequalities.  
For rectangular
sections, the position constraints are
\begin{align}
  \Big( \frac{2 \bar{w}_i(s) - (\bar{w}_{i,max}(s)+\bar{w}_{i,min}(s))}{(\bar{w}_{i,max}(s)-\bar{w}_{i,min}(s))} \Big)^2 -1 \leq 0,
\label{eq:rect_constr} 
\end{align}
$\forall s \in [0, L]$, $\forall i = 1,2,$ where the non constant functions $\bar{w}_{i,min}(\cdot)$ and
$\bar{w}_{i,max}(\cdot)$ identify the obstacle boundaries.  For circular
sections, the position constraint is
\begin{align}
  \Big( \frac{||\bar{\pmb{w}}(s)||}{\bar{r}_{obs}(s)} \Big)^2 -1 \leq 0, \quad \forall s \in [0, L],
  \label{eq:circ_constr}
\end{align}
where the non constant function $\bar{r}_{obs}(\cdot)$ identifies the circular
boundary of the collision free region. 

The new formulation \eqref{eq:mintime} is equivalent to \eqref{eq:mintime_standard} since
trajectories solving \eqref{eq:mintime} can be mapped into trajectories solving
\eqref{eq:mintime_standard}.  Using transverse
coordinates, the formulation of position constraints is simpler than
using the time-parameterized inertial frame position. Since collision free region boundaries $\bar{w}_{i,min}(\cdot)$,
$\bar{w}_{i,max}(\cdot)$ and $\bar{r}_{obs}(\cdot)$ are arc-length functions,
even complex spaces can be modelled. 
\subsection{Numerical solution to the optimization problem}
In order to solve problem \eqref{eq:mintime}, we use a combination of the
PRojection Operator based Newton method for Trajectory Optimization (PRONTO)
\cite{hauser2002projection} with a barrier function approach
\cite{hauser2006barrier}.  We relax state-input constraints by adding them in
the cost functional, i.e., given the state-input curve
$\xi = (\bar{\pmb{q}}_S(\cdot),\bar{\pmb{u}}(\cdot))$, we consider the problem
\begin{equation}
  \begin{split}
    \min_{\xi} &\;\;  h(\xi; \bar{k}) + \epsilon b_\nu(\xi)\\
    \!\!\subj &\;\; \bar{\pmb{q}}'_S(s) =
    \bar{f}(\bar{\boldsymbol{q}}_S(s),\bar{\boldsymbol{u}}(s)), \quad
    \bar{\pmb{q}}_S(0) = \pmb{q}_{S0}
  \end{split}
  \label{eq:mintime2}
\end{equation}
where
\begin{align}
  h(\xi; \bar{k})&=\int_0^L \!\!\! \quad \frac{1-\bar{k}(s) \bar{w}_1(s)}{\bar{\text{v}}_t(s)} \; ds \nonumber\\
  b_\nu(\xi) &= 
               \int_0^L \sum_j \beta_\nu (-c_j(\bar{\pmb{q}}_S(s),\bar{\pmb{u}}(s))) \; ds \nonumber\\
  \beta_\nu (x) &=  
                  \begin{cases} 
                    -\log(x) & x > \nu\\
                    -\log(\nu) + \frac{1}{2} \big[ (\frac{x - 2\nu}{\nu})^2 -1
                    \big] & x \leq \nu,
                  \end{cases} \nonumber
\end{align}
and $c_j(\bar{\pmb{q}}_S(s),\bar{\pmb{u}}(s))$, j=1,..,7 denote the
constraints in problem \eqref{eq:mintime}.  The strategy to find a solution to
\eqref{eq:mintime} consists in iteratively solving problem \eqref{eq:mintime2}
by reducing the parameters $\epsilon$ and $\nu$ at each iteration, and thus pushing
the trajectory towards the constraint boundaries.
Each instance of problem \eqref{eq:mintime2} is solved by means of the PRONTO
algorithm.  This is based on a properly designed \emph{projection operator}
$\mathcal{P} : \xi_c \rightarrow \xi$, mapping a state-control curve
$\xi_c =(\bar{\pmb{q}}_c(\cdot), \bar{\pmb{u}}_c(\cdot))$ into a system
trajectory $\xi$, by the nonlinear feedback system
\begin{align}
  \bar{\pmb{q}}'_S(s) & = \bar{f}(\bar{\boldsymbol{q}}_S(s),\bar{\boldsymbol{u}}(s)), \quad \bar{\pmb{q}}_S(0) = \pmb{q}_{S0} \nonumber\\
  \bar{\pmb{u}}(s) & = \bar{\pmb{u}}_c(s) + \bar{K}(s)(\bar{\pmb{q}}_c(s)-\bar{\pmb{q}}_S(s)),
                     \label{eq:proj_oper_def}
\end{align}
where the feedback gain $\bar{K}(\cdot)$ is designed by solving a suitable
linear quadratic optimal control problem.  The projection operator is used to
convert the dynamically constrained optimization problem \eqref{eq:mintime2}
into the unconstrained problem
\begin{equation}
  \begin{split}
    \min_{\xi} &\; g(\xi;\bar{k}),
  \end{split}
  \label{eq:mintime3}
\end{equation}
where
$g(\xi;\bar{k}) = h(\mathcal{P}(\xi); \bar{k}) + \epsilon
b_\nu(\mathcal{P}(\xi)).$
Then, using an (infinite dimensional) Newton descent method, a local minimizer
of \eqref{eq:mintime3} is computed.  

An initial trajectory, with a position matching the ``reference" path, is used
to initialize the algorithm.  Given the current trajectory iterate $\xi_i$, the search
direction $\zeta_i$ is obtained by solving a linear quadratic optimal control
problem with cost
$Dg(\xi_i; \bar{k}) \cdot \zeta + \frac{1}{2} D^2 g(\xi_i;
\bar{k})(\zeta,\zeta)$,
where $\zeta \mapsto Dg(\xi_i; \bar{k}) \cdot \zeta$ and
$\zeta \mapsto D^2 g(\xi_i; \bar{k})(\zeta,\zeta)$ are respectively the first
and second Fr\'echet differentials of the functional $g(\xi,\bar{k})$ at
$\xi_i$. Then, the curve $\xi_i + \gamma_i \zeta_i$, where $\gamma_i$ is a step
size obtained through a standard backtracking line search, is projected, by
means of the projection operator, in order to get a new trajectory $\xi_{i+1}$.
The strength of this approach is that the local minimizer of \eqref{eq:mintime3}
is obtained as the limit of a sequence of trajectories, i.e., curves satisfying
the dynamics.  Furthermore, the feedback system \eqref{eq:proj_oper_def},
defining the projection operator, allows us to generate trajectories in a
numerically stable manner.

\section{Numerical computations}
\label{sec:computations}
In this section, we present numerical computations to show the effectiveness of
the proposed strategy.  

We consider two different scenarios with
a nano-quadrotor having mass
$m = 0.021 \; \text{kg}$.

\begin{figure}[htbp]
  \begin{center}
    \subfloat[][Path]
    {\includegraphics[scale=0.35]{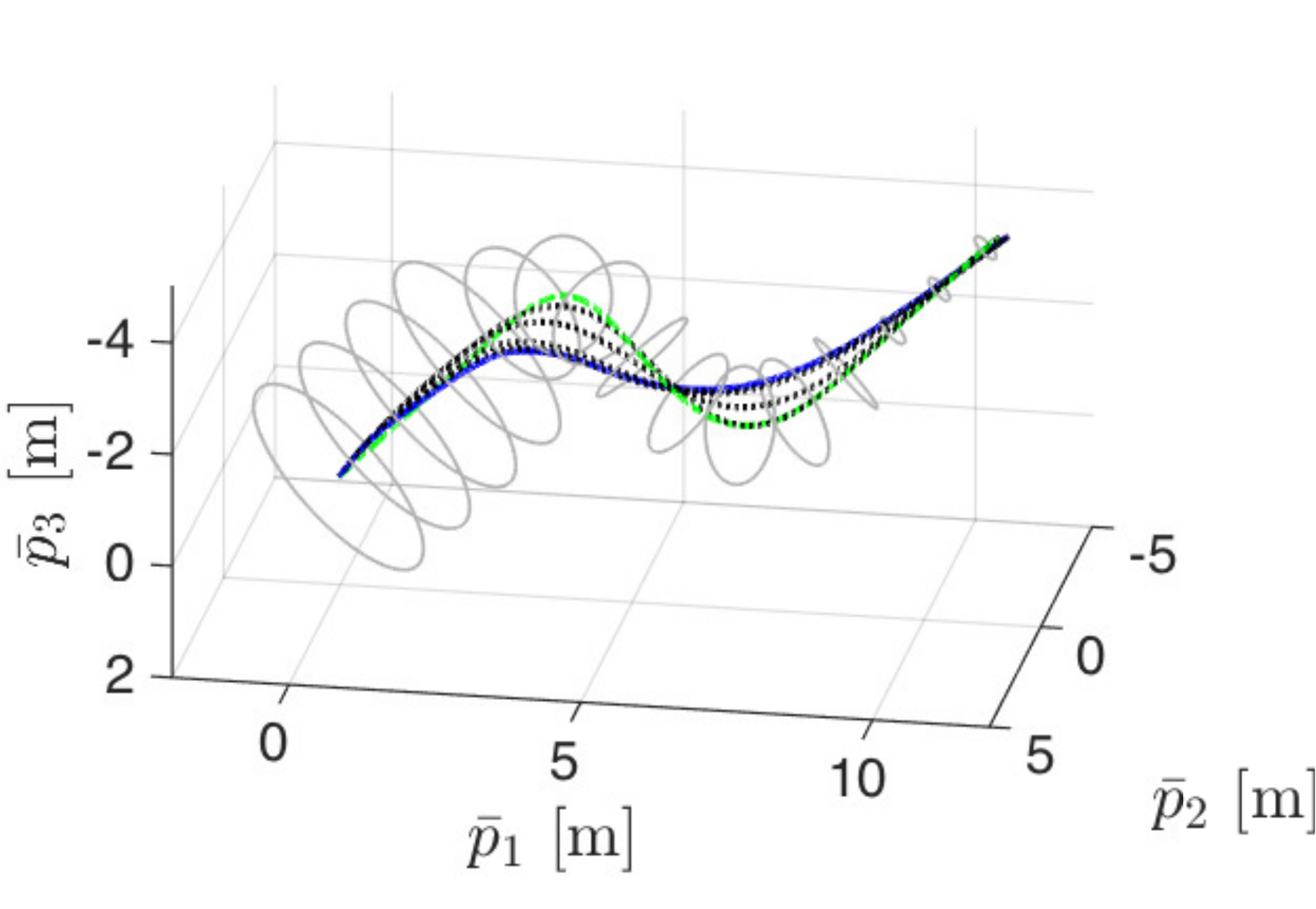}\vspace{0.5cm}\label{fig:path_s1}}\\
    \subfloat[][Distance from the ``reference" path]
    {\hspace{-0.1cm}\includegraphics[width=5.0cm]{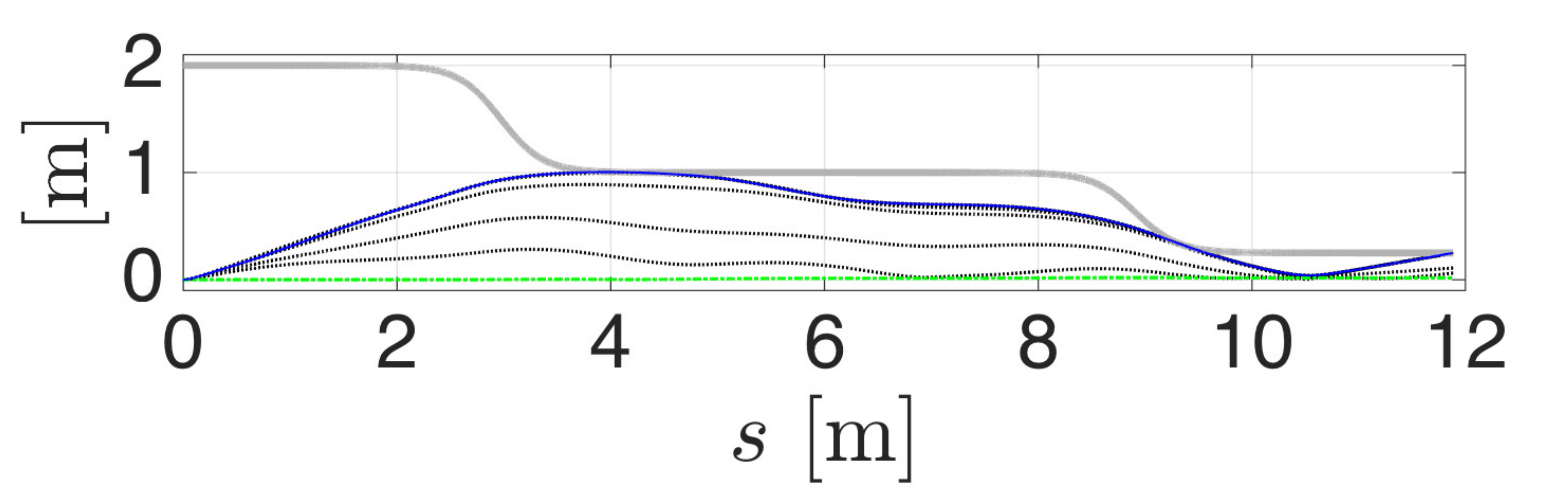}\label{fig:dd_s1}}\\
        \subfloat[][Velocity $\bar{\text{v}}_t$]
    {\hspace{-0.1cm}\includegraphics[width=3.9cm]{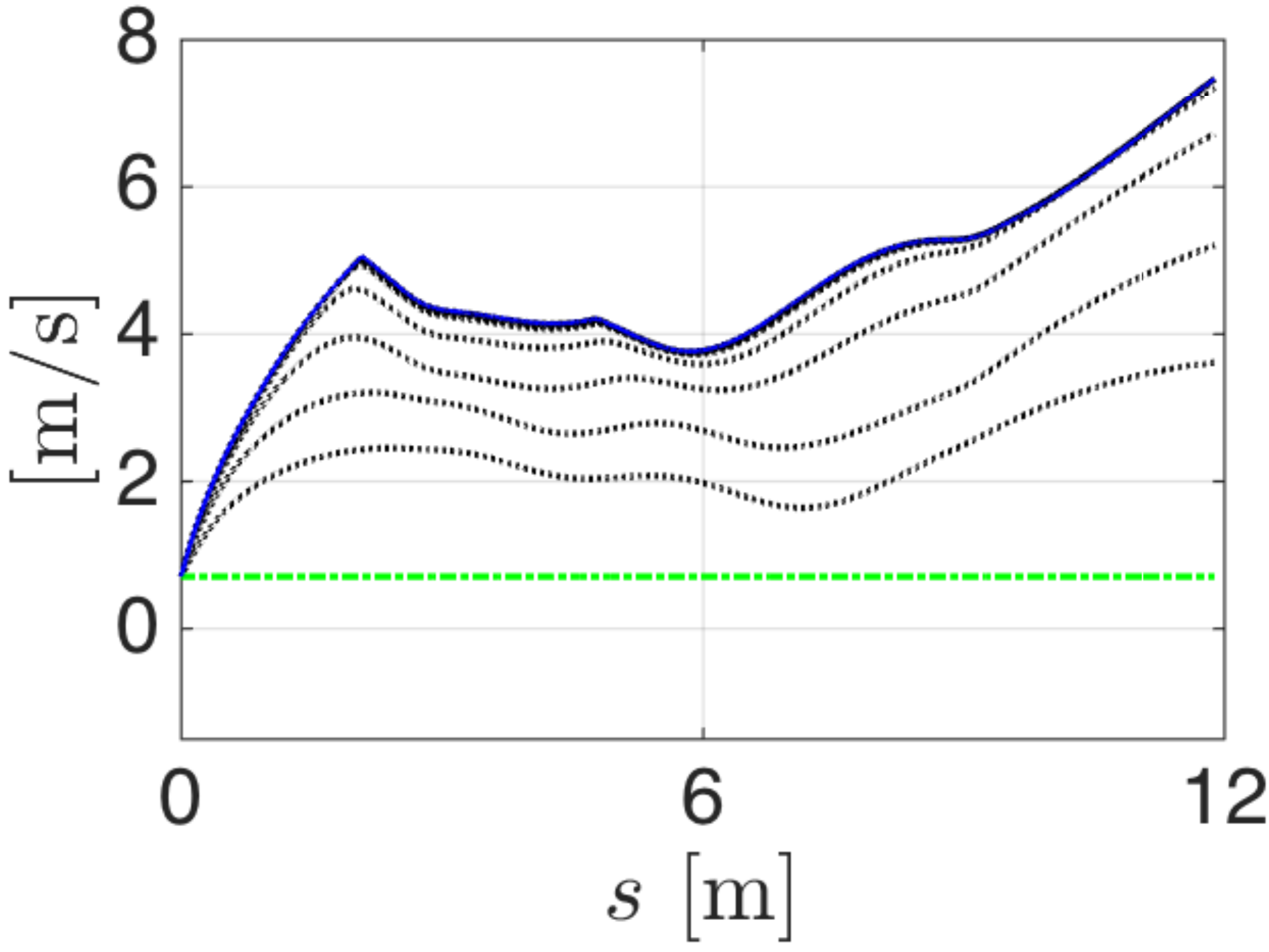}\label{fig:dp1_s1}}
    \subfloat[][Roll angle $\bar{\varphi}$]
    {\hspace{0.36cm}\includegraphics[width=3.9cm]{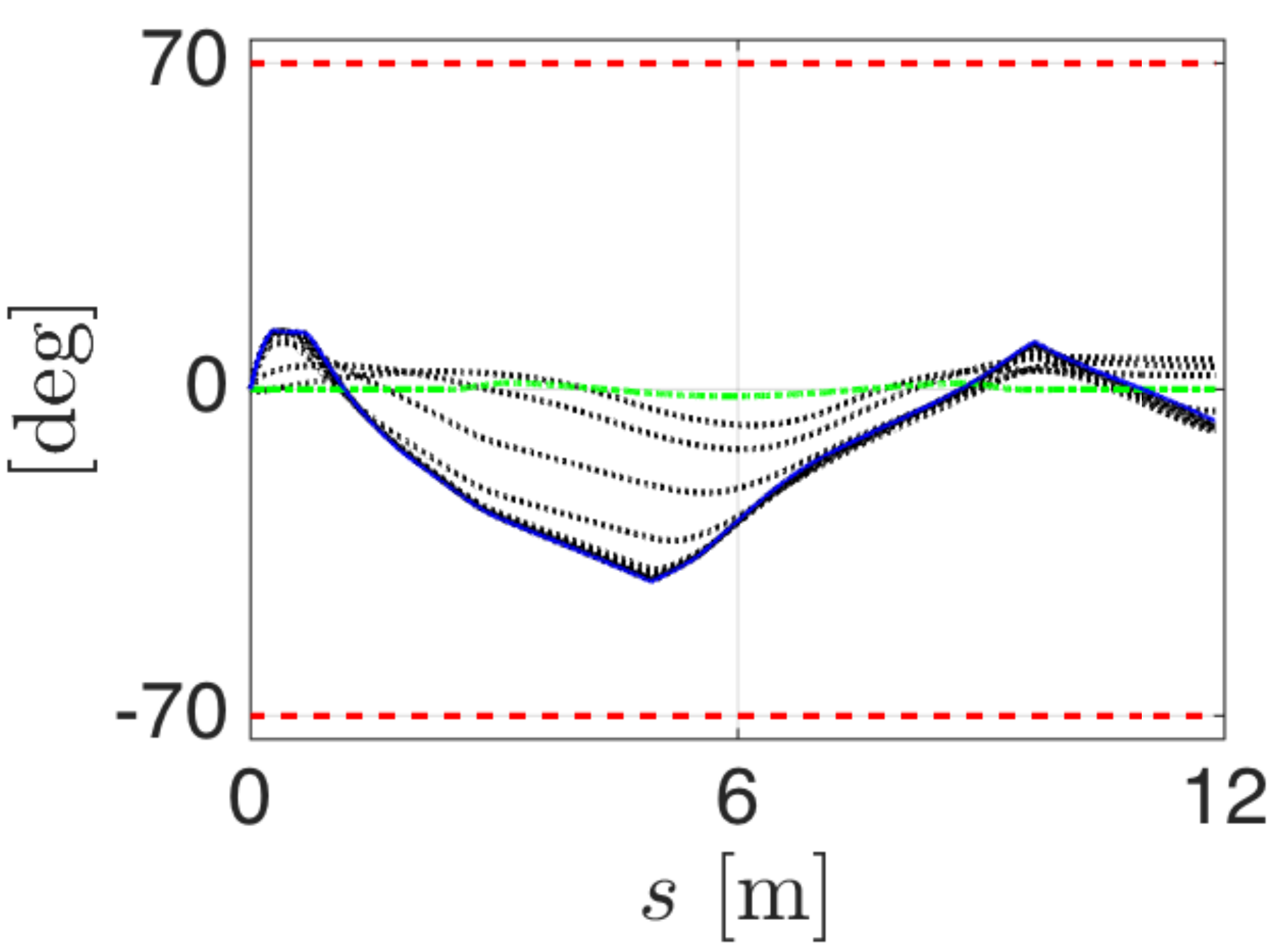}\label{fig:phi_s1}}\\
    \caption{Scenario $1$ (path, distance from the ``reference" path and states). Initial
      (dot-dashed green), intermediate (dotted black) and minimum-time (solid
      blue) trajectory. Tube boundaries (gray) and angular constraints (dashed red).}
    \label{fig:scenario1_path}
  \end{center}
\end{figure}

\begin{figure}[ht!]
  \begin{center}
    \subfloat[][Roll rate $\bar{\omega}_1$]
    {\hspace{-0.1cm}\includegraphics[width=3.9cm]{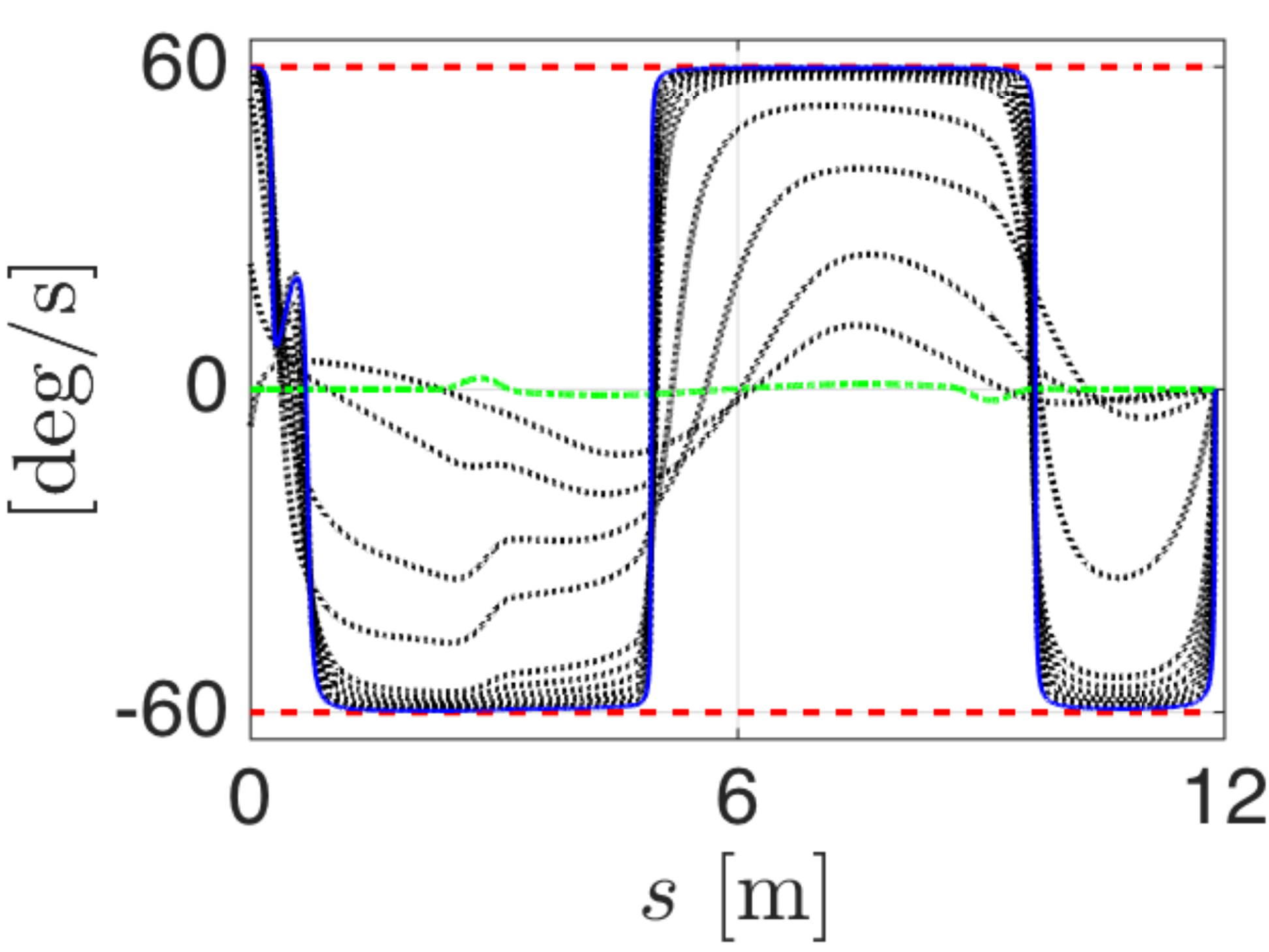}\label{fig:pp_s1}}
    \subfloat[][Pitch rate $\bar{\omega}_2$]
    {\hspace{0.36cm}\includegraphics[width=3.9cm]{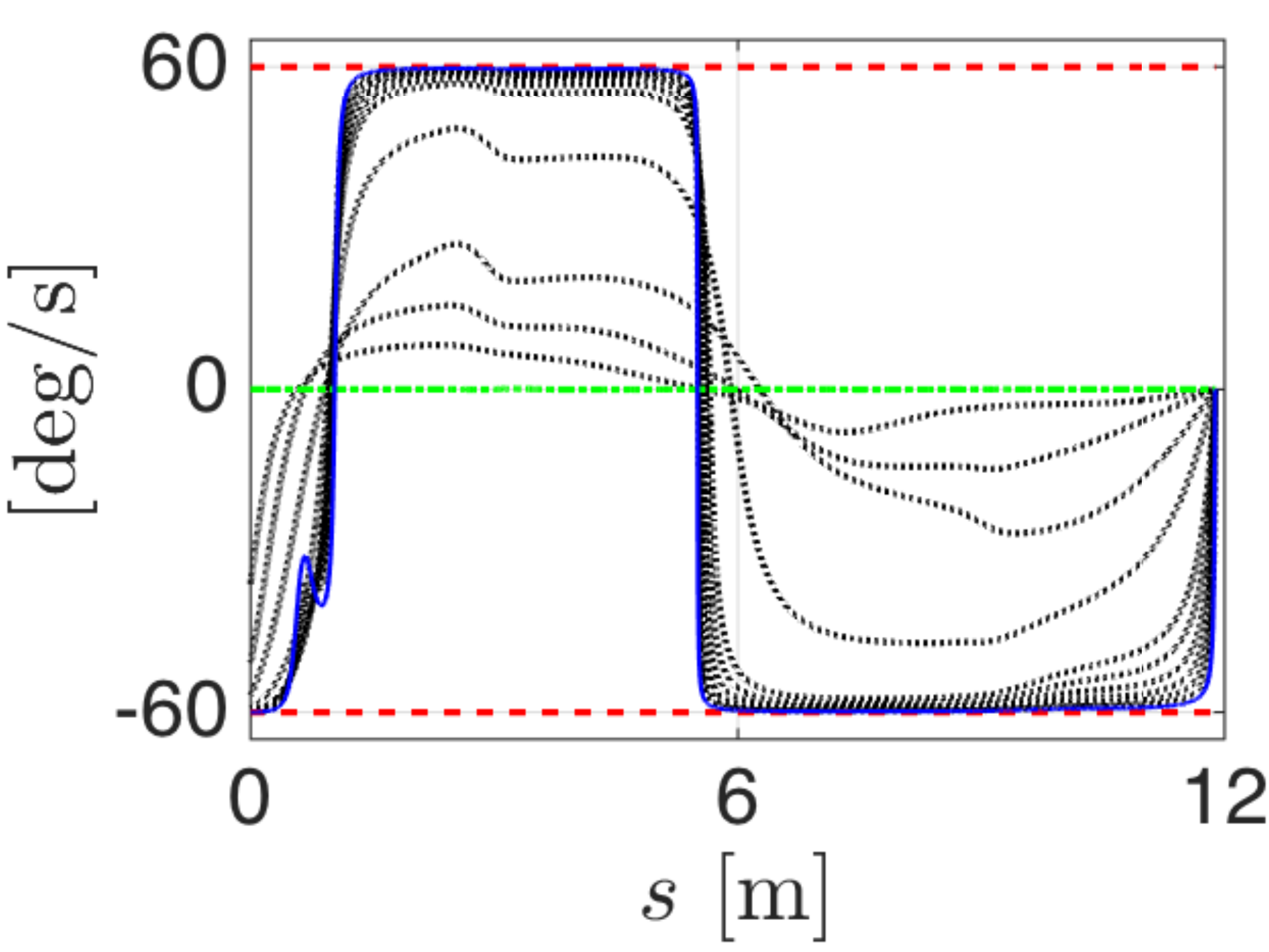}\label{fig:qq_s1}}\\
    \subfloat[][Yaw rate $\bar{\omega}_3$]
    {\hspace{-0.1cm}\includegraphics[width=3.9cm]{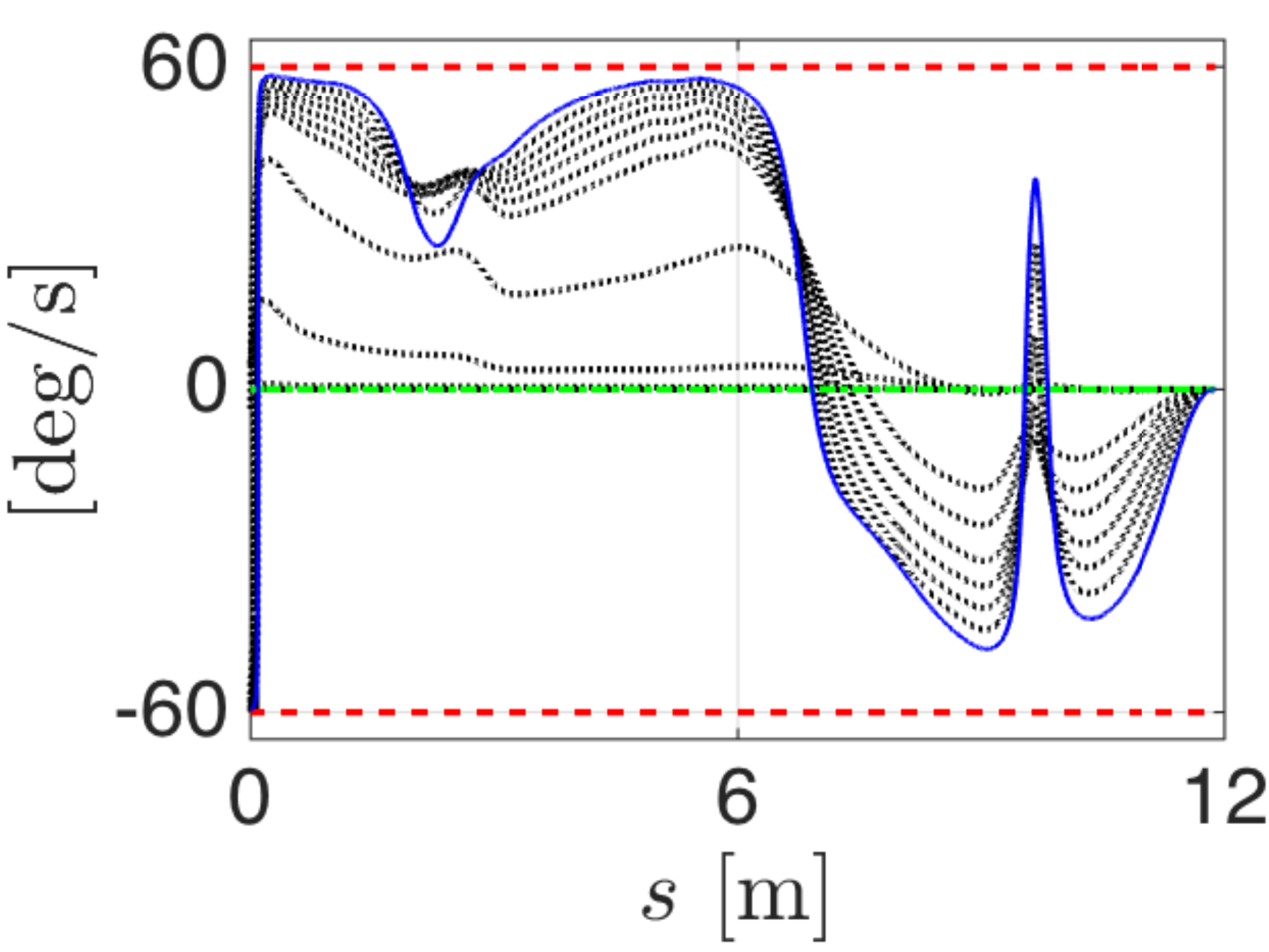}\label{fig:rr_s1}}
    \subfloat[][Thrust $\bar{f}$]
    {\hspace{0.36cm}\includegraphics[width=3.9cm]{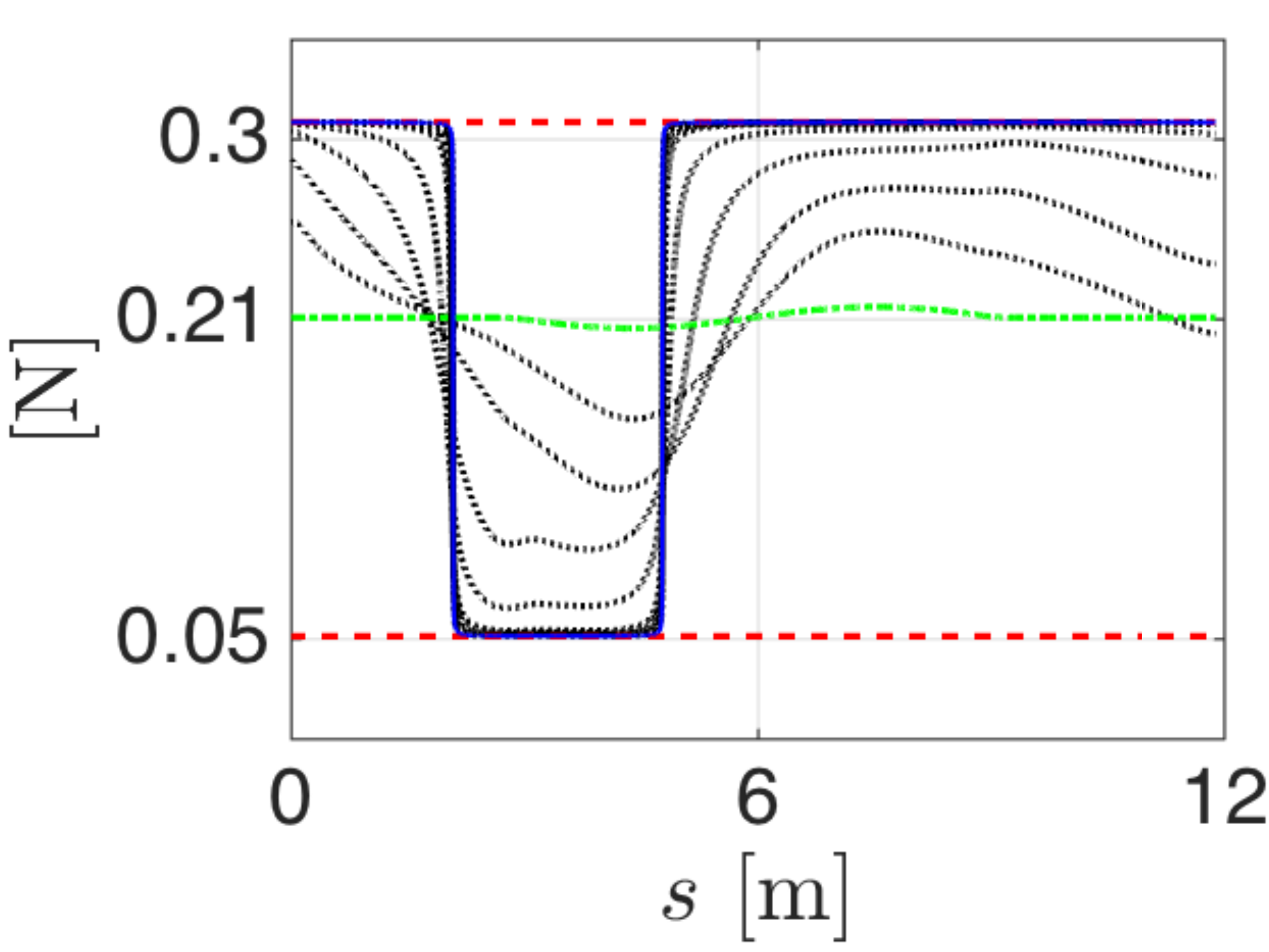}\label{fig:ff_s1}}\\
    \caption{Scenario $1$ (inputs). Initial (dot-dashed green), intermediate
      (dotted black) and minimum-time (solid blue) trajectory. Constraint
      boundaries (dashed red).}
    \label{fig:scenario1_inputs}
  \end{center}
\end{figure}

In the first scenario, the vehicle moves inside a tube with a $3$D center line and circular
transversal sections getting smaller along it.
Results are depicted in Figures \ref{fig:scenario1_path} and
\ref{fig:scenario1_inputs}.  
As regards the initial trajectory (dot-dashed green line), by choosing as ``reference" path the 3D center-line of the tube, the collision free
region is defined by constraint \eqref{eq:circ_constr} where
$\bar{r}_{obs}(\cdot)$ is depicted in Figure \ref{fig:dd_s1} (gray line).
Furthermore, we set an initial constant velocity 
and a zero yaw angle. The remaining
initial states and inputs are computed by tracking these positions and
velocities. 
Note that, in real applications, the bounds of the collision free region should be restricted (with respect to the physical boundaries of the region) in order to consider the vehicle dimension and an imperfect tracking of the trajectory.
Each intermediate optimal trajectory (dotted black line) is computed by iteratively solving problem
\eqref{eq:mintime2} with decreasing values of $\epsilon$
and $\nu$.  
As regards the minimum-time trajectory (blue line), the maneuver is
performed in $0.36 \; \text{s}$, the optimal path touches the tube at arc-lengths $4.11 \;\text{m}$, $9.35 \;\text{m}$ and $11.88 \;\text{m}$ (Figure
\ref{fig:dd_s1}) and the velocity $\bar{\text{v}}_t$ (Figure \ref{fig:dp1_s1})
reaches more than $7 \; \text{m/s}$ at the end of the path.
As regards inputs, while constraints on thrust
(Figure \ref{fig:ff_s1}) are always active, roll and pitch rate (Figures
\ref{fig:pp_s1} and \ref{fig:qq_s1}, respectively) are saturated after
$s = 1 \; \text{m}$ and the yaw rate (Figure \ref{fig:rr_s1}) touches the
constraint boundary only in the beginning of the maneuver.
%
	
\begin{figure}[htbp]
		\begin{center}
			\vspace{-0.2cm}
			\subfloat[][Path: 3D view]
			{\includegraphics[scale=0.3]{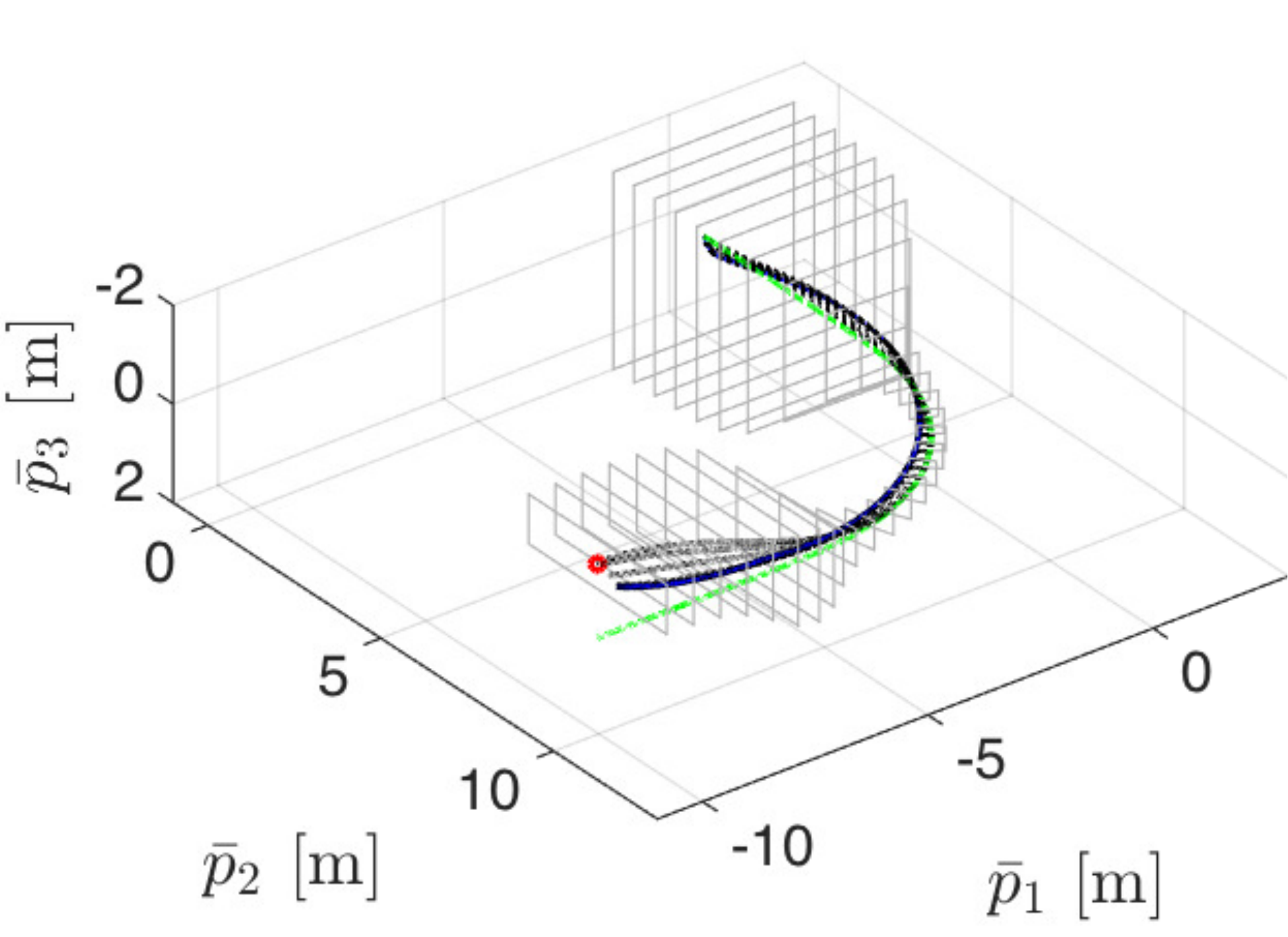}\vspace{0.25cm}\label{fig:path_s2}}\\
			\subfloat[][Transverse coordinate $\bar{w}_1$]
			{\hspace{-0.1cm}\includegraphics[width=3.9cm]{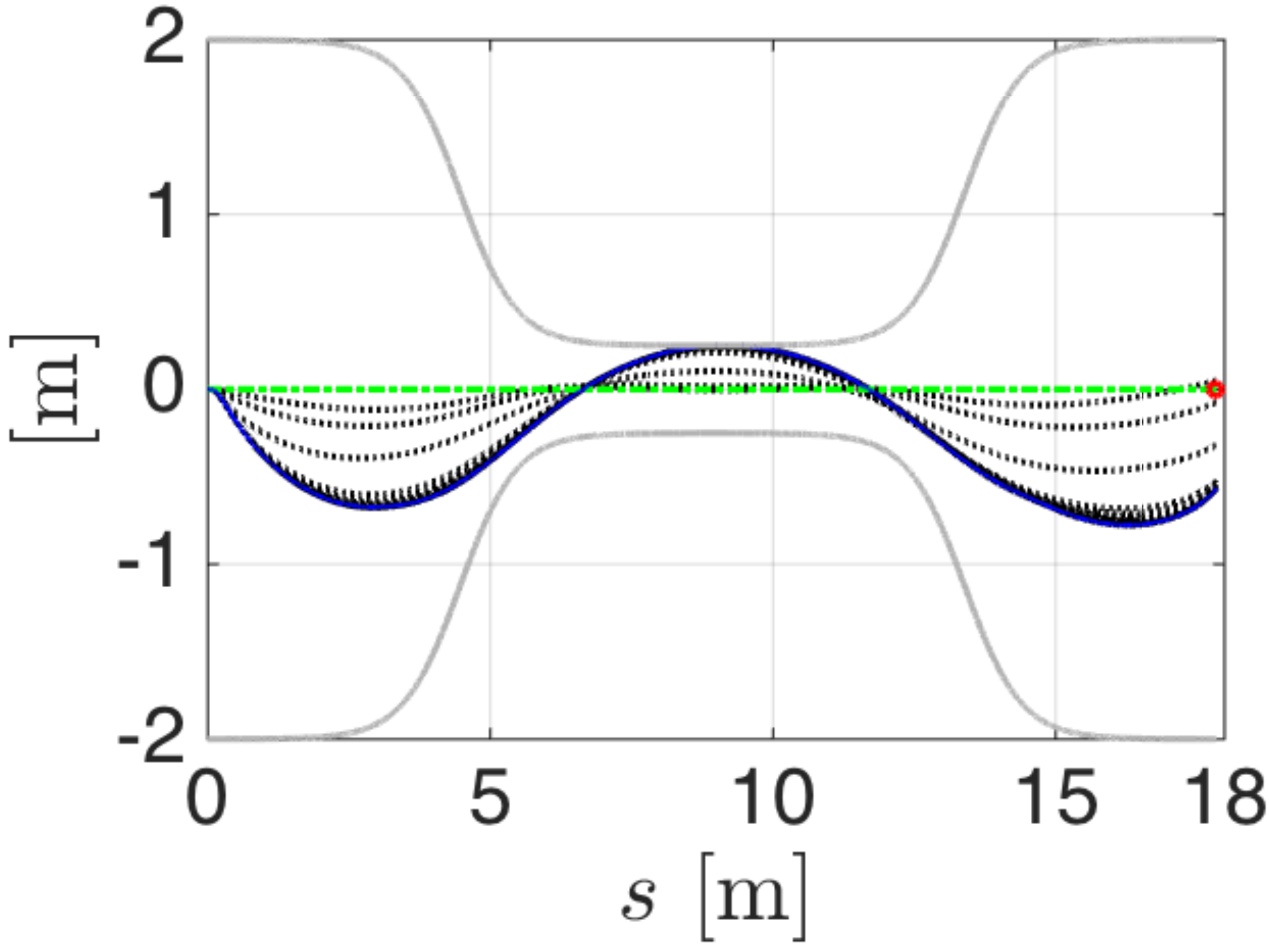}\label{fig:w1_s2}}
			\subfloat[][Transverse coordinate $\bar{w}_2$]
			{\hspace{0.36cm}\includegraphics[width=3.9cm]{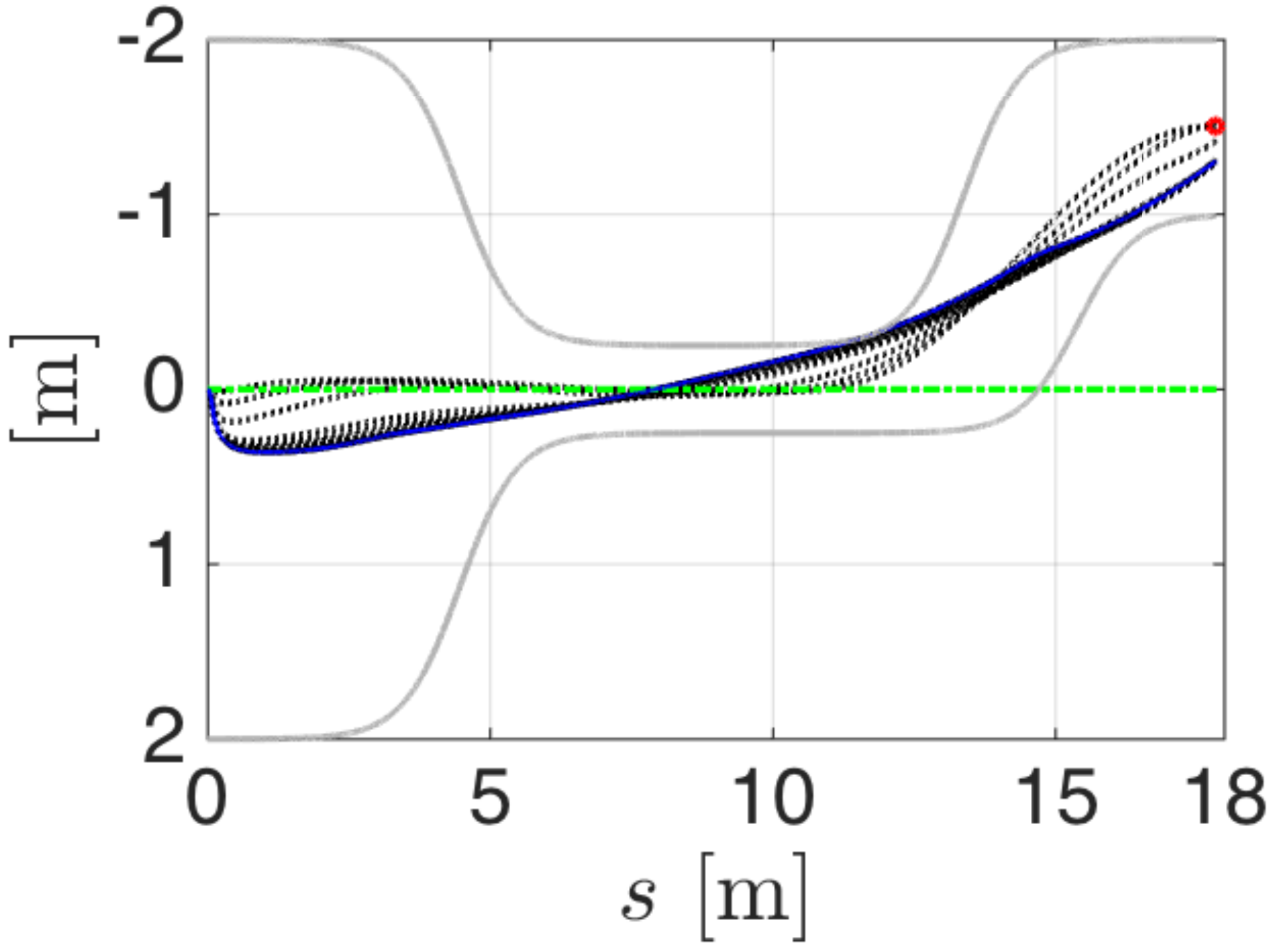}\label{fig:w2_s2}}\\
						\subfloat[][Velocity $\bar{\text{v}}_t$]
			{\hspace{-0.1cm}\includegraphics[width=3.9cm]{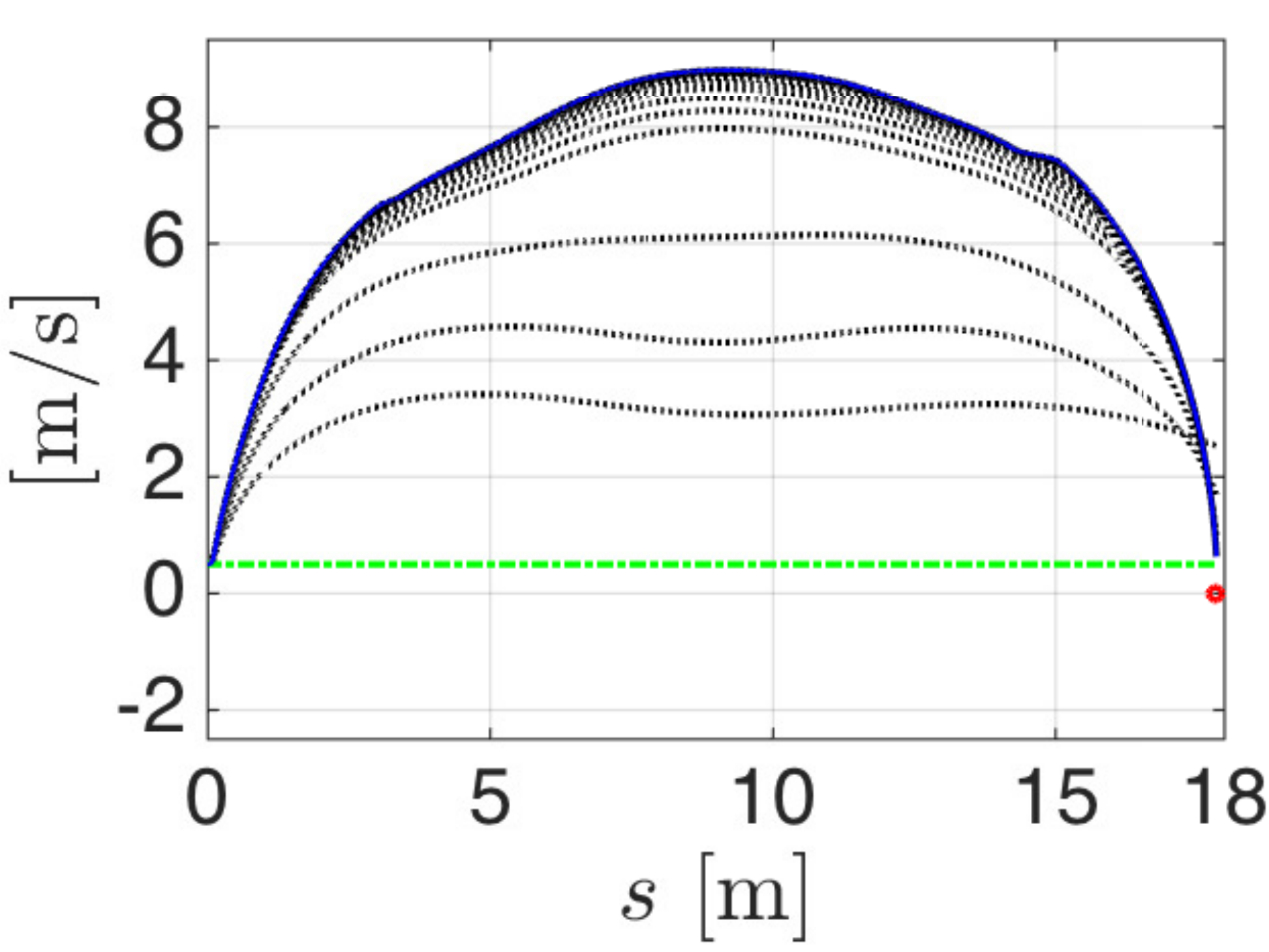}\label{fig:dp1_s2}}
			\subfloat[][Roll angle $\bar{\varphi}$]
			{\hspace{0.36cm}\includegraphics[width=3.9cm]{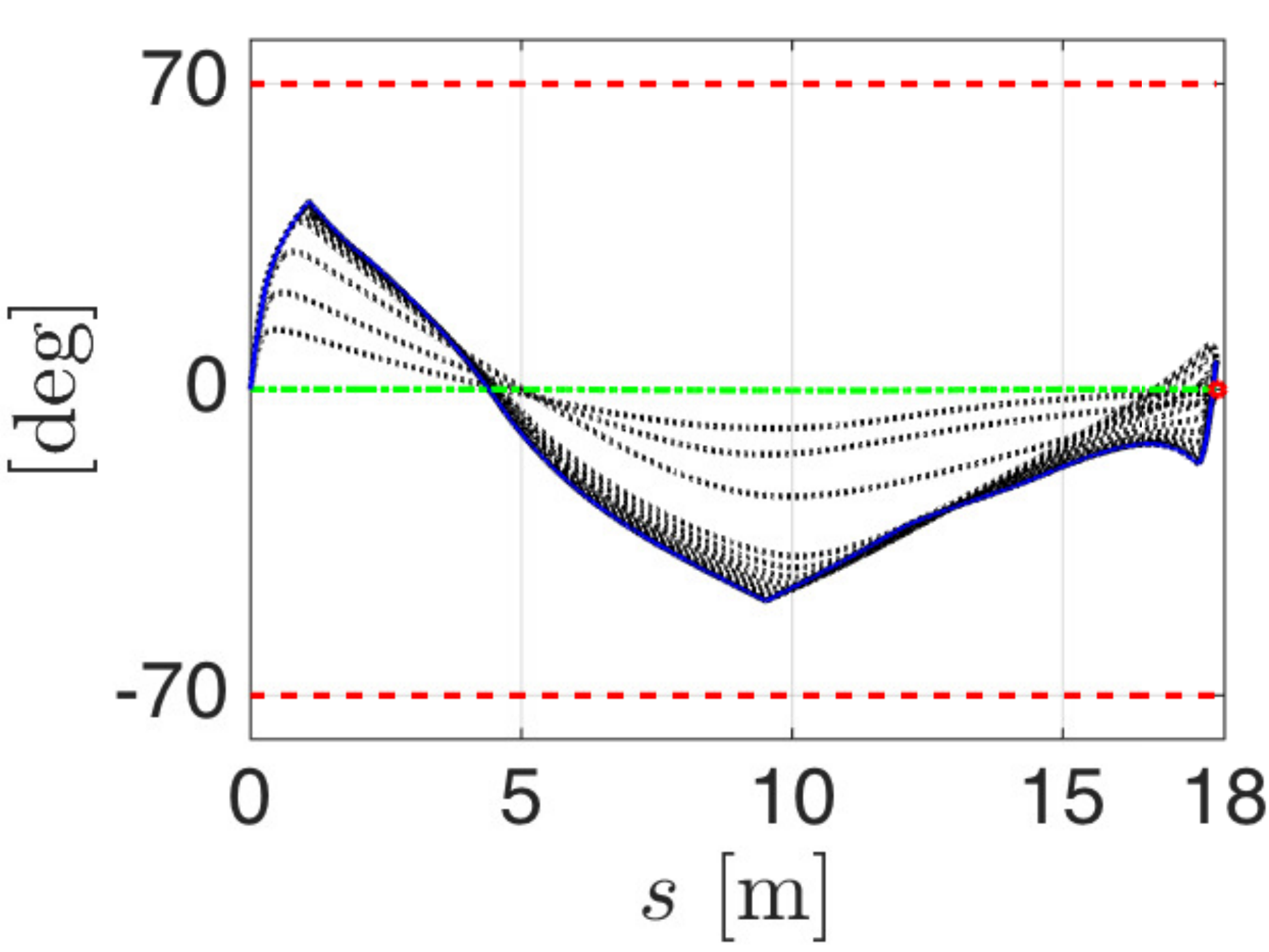}\label{fig:phi_s2}}\\
			\caption{Scenario $2$ (path and states). Initial (dot-dashed green), intermediate (dotted black) and minimum-time (solid blue) trajectory. Collision free region boundary (gray), angular constraints (dashed red) and final desired state (red circle).}
			\label{fig:scenario2_path}
		\end{center}
	\end{figure}

\begin{figure}[htbp]
		\begin{center}
			\subfloat[][Roll rate $\bar{\omega}_1$]
			{\hspace{-0.1cm}\includegraphics[width=3.9cm]{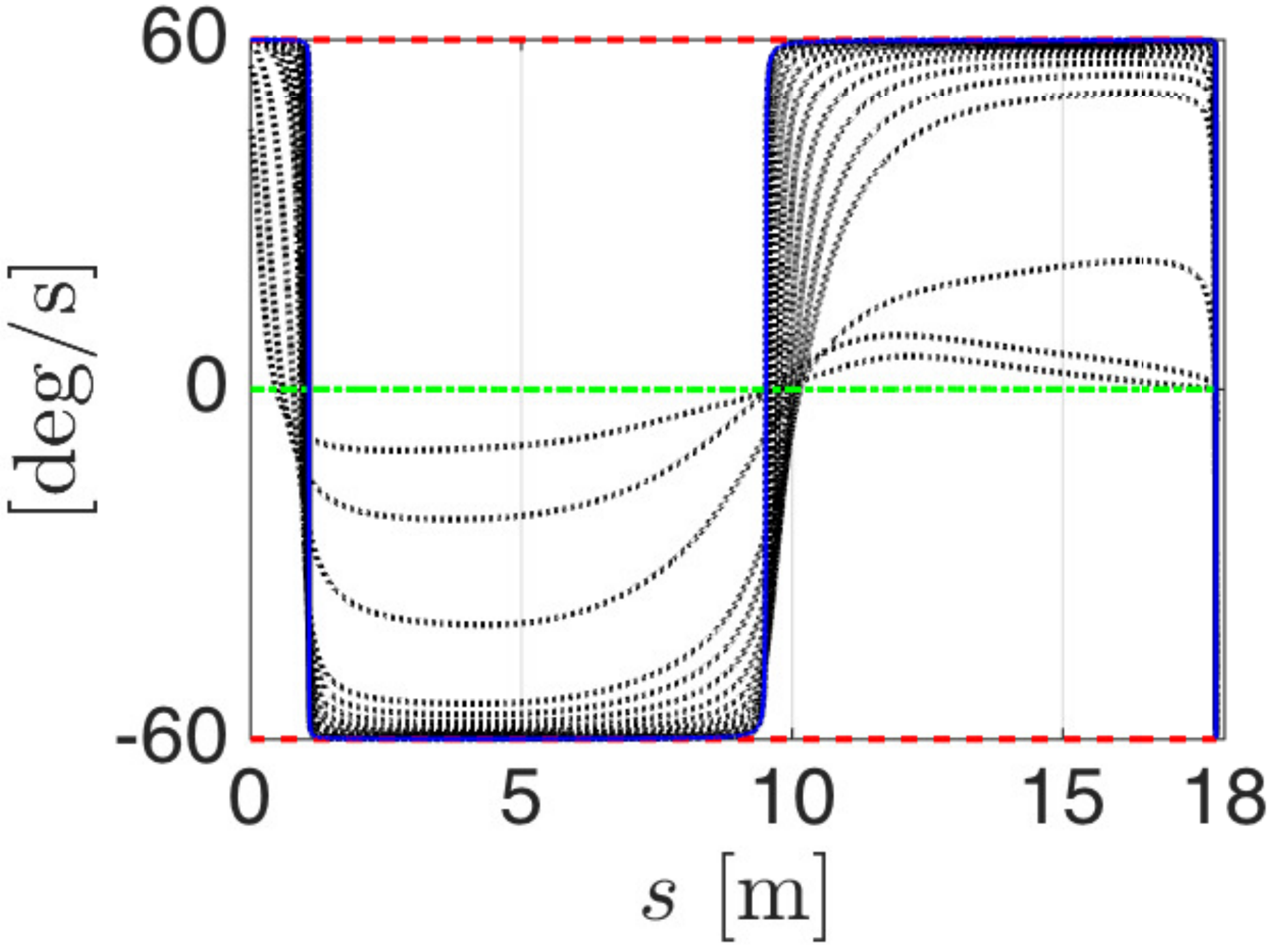}\label{fig:pp_s2}}
			\subfloat[][Pitch rate $\bar{\omega}_2$]
			{\hspace{0.36cm}\includegraphics[width=3.9cm]{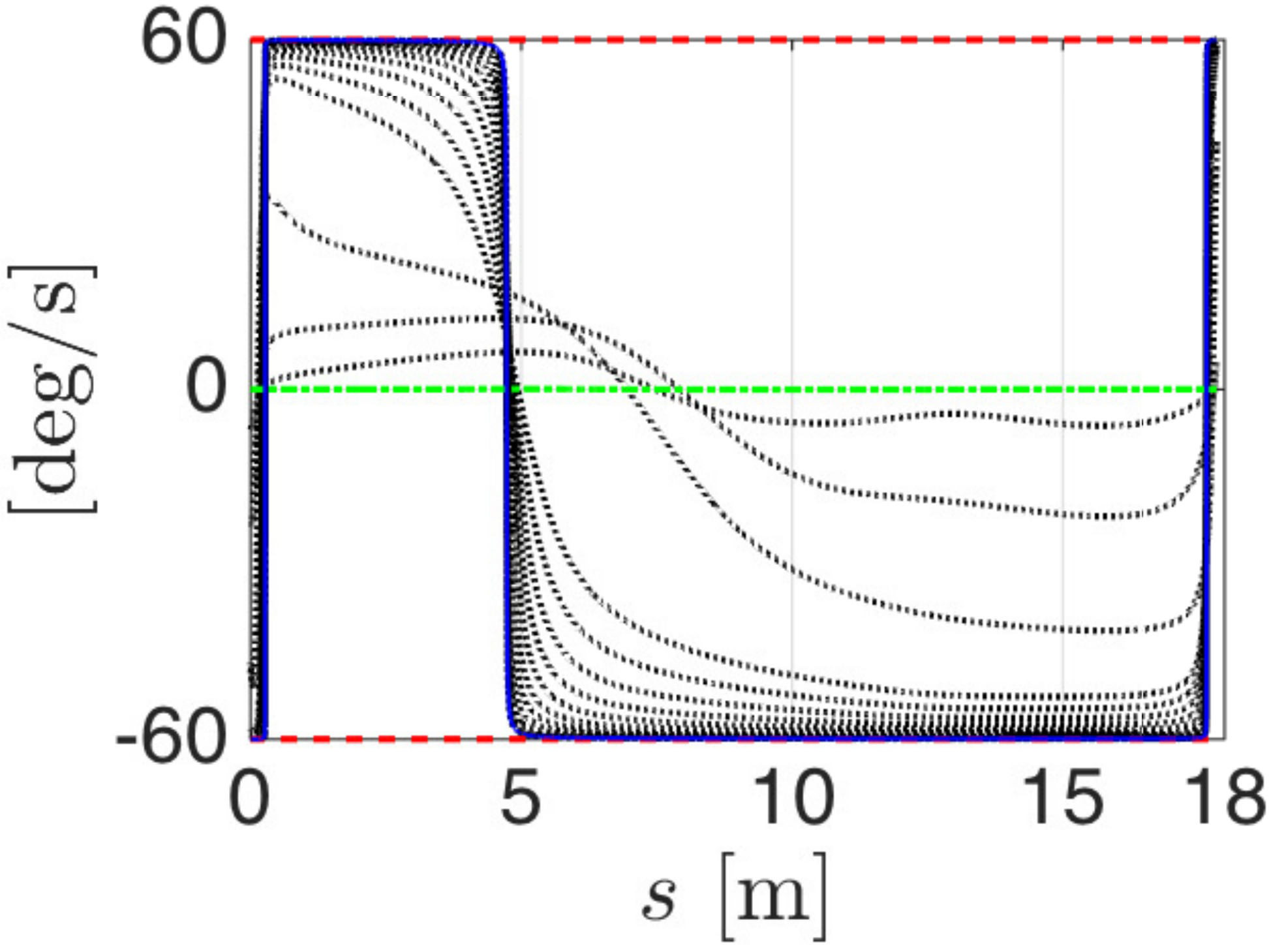}\label{fig:qq_s2}}\\
			\subfloat[][Yaw rate $\bar{\omega}_3$]
			{\hspace{-0.1cm}\includegraphics[width=3.9cm]{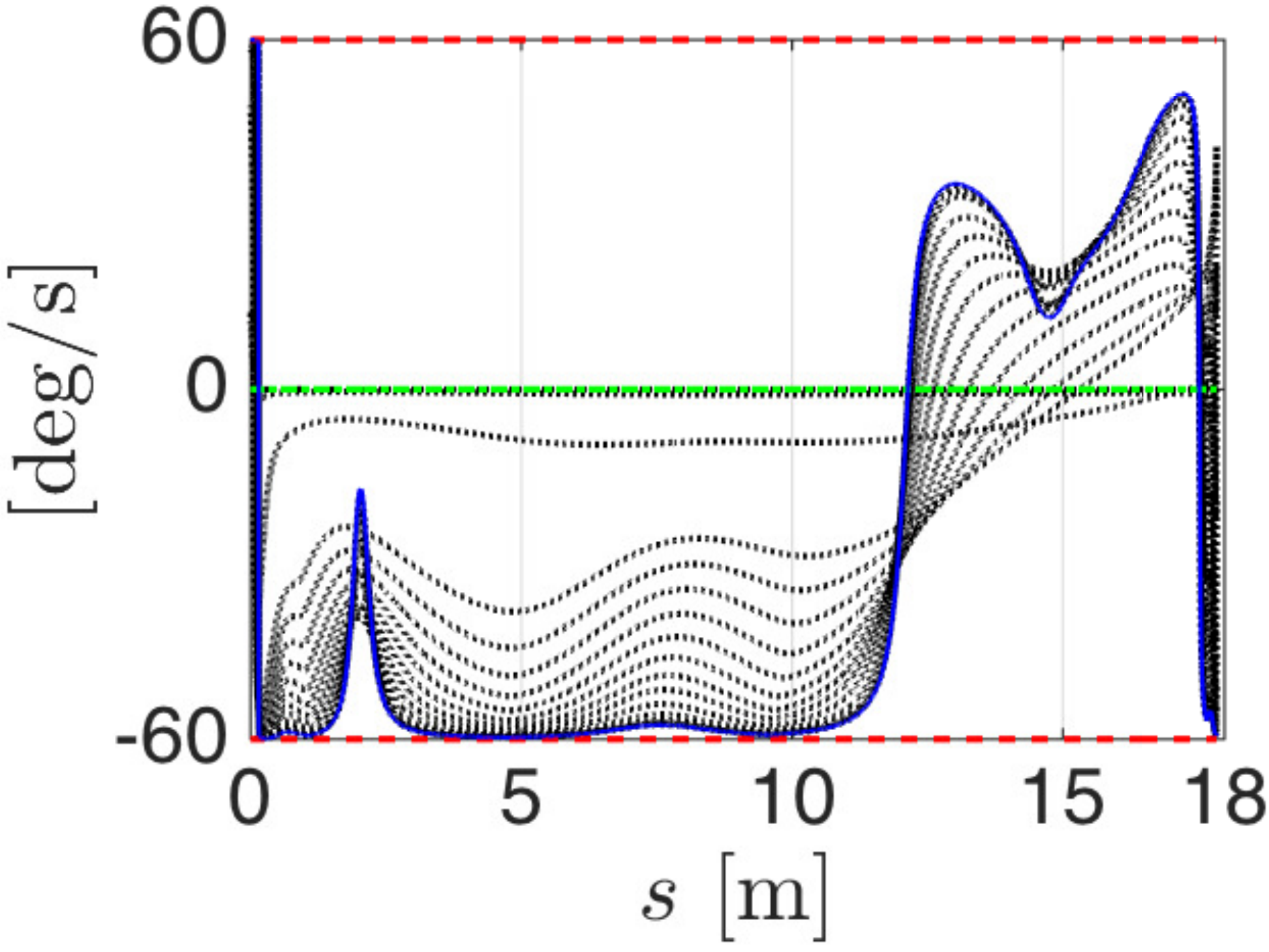}\label{fig:rr_s2}}
			\subfloat[][Thrust $\bar{f}$]
			{\hspace{0.36cm}\includegraphics[width=3.9cm]{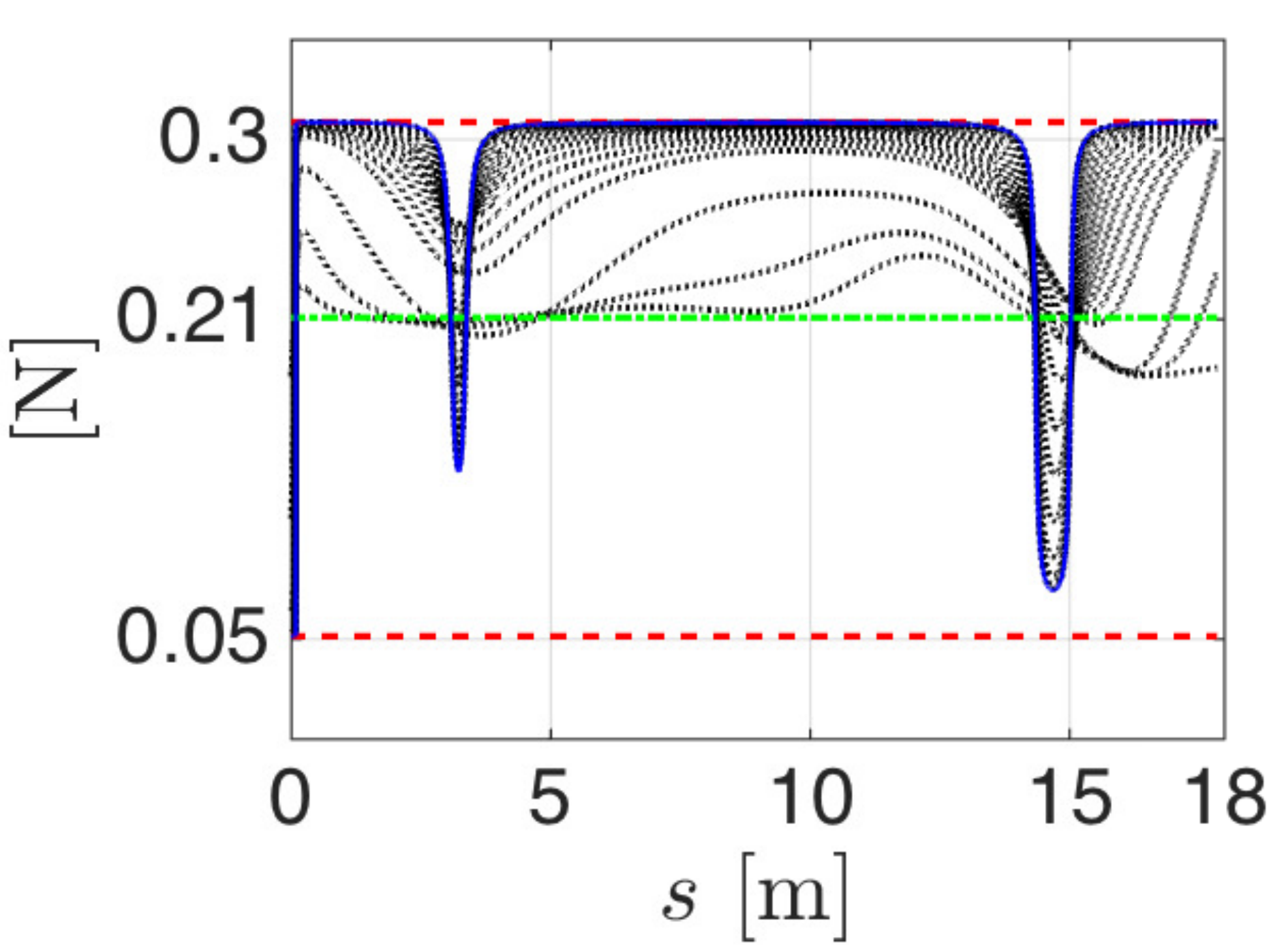}\label{fig:ff_s2}}\\						
			\caption{Scenario $2$ (inputs). Initial (dot-dashed green), intermediate (dotted black) and minimum-time (solid blue) trajectory. Constraint boundaries (dashed red).}
			\label{fig:scenario2_inputs}
		\end{center}
	\end{figure}

In the second scenario, the vehicle moves from one room to another through a
narrow corridor with a squared section with side length $0.5 \; \text{m}$.
Leaving the corridor, the vehicle encounters an obstacle and is forced to
fly above it, having a final altitude $\bar{p}_3(L) \leq -1 \; \text{m}$. 
In this scenario, as an additional
requirement, we seek for a trajectory reaching the desired final state
$\pmb{\bar{q}}_{S}^d = [0 \; -1.5 \; 0 \; 0 \; 0 \; 0 \; 0 \; 0]^T$.  In order
to fulfil this objective, we use a penalty function approach by adding the term
$\frac{1}{2} (\pmb{\bar{q}}_{S}^d-\pmb{\bar{q}}_{S}(L))^T \rho I \;
(\pmb{\bar{q}}_{S}^d-\pmb{\bar{q}}_{S}(L))$,
where $\rho$ is a penalty parameter and $I$ is the identity $8 \times 8$ matrix,
in the cost functional of problem \eqref{eq:mintime}.  Results are depicted in
Figures \ref{fig:scenario2_path} 
and
\ref{fig:scenario2_inputs}.
We choose as ``reference" path an unfeasible curve on the $\bar{p}_1 - \bar{p}_2$
plane and
the collision free region is defined by constraint \eqref{eq:rect_constr} where
obstacle boundaries $\bar{w}_{i,min}(\cdot)$ and $\bar{w}_{i,max}(\cdot)$,
$i = 1,2$ are depicted in gray in Figures \ref{fig:w1_s2} and \ref{fig:w2_s2}.
In particular, $\bar{w}_{2,min}$ is defined in order to take into account the presence of the obstacle that forces the vehicle to fly above it.
 An initial constant velocity
and a zero yaw angle are considered and the remaining initial states and inputs are computed as in the
first scenario.
Note that, even if the initial trajectory is unfeasible, 
intermediate trajectories are feasible.  As regards the minimum-time
trajectory, the maneuver is performed in $0.5 \; \text{s}$ and the path touches
the constraint boundaries when the vehicle is inside the corridor (Figures \ref{fig:w1_s2} and \ref{fig:w2_s2}). The velocity $\bar{\text{v}}_t$ (Figure
\ref{fig:dp1_s2}) reaches a peak of $9 \; \text{m/s}$ in the middle of the path
and approaches zero at the end. 
As regards
the inputs, only constraints on roll and pitch rates (Figures \ref{fig:pp_s2},
\ref{fig:qq_s2}, respectively) are always active.
Furthermore, note that the desired final state
is not exactly reached, since our penalty function approach only
guarantees approximated solutions. An ad-hoc strategy should be designed to
exactly reach the desired final state.
	%
\section{Conclusion}
In this paper, we have presented a strategy to compute minimum-time trajectories
for quadrotors in constrained spaces. Our approach consists of: (i)
generating a ``reference" path, (ii) expressing the vehicle dynamics
using transverse coordinates, and (iii)
re-defining constraints in the new coordinates. We obtain a re-formulation
of the problem, which we solve by combining the PRONTO algorithm with a barrier
function approach. Numerical computations on two challenging scenarios prove the
effectiveness of the strategy.
As a future work, we aim to compare our numerical
results with a theoretical analysis on particular scenarios.
\bibliographystyle{IEEEtran}
\bibliography{bibliography}  
\end{document}